\input amstex
\magnification=1200
\documentstyle{amsppt}
\NoBlackBoxes
\NoRunningHeads
\topmatter
\title Interactive games, dialogues and the verbalization
\endtitle
\author Denis V. Juriev
\endauthor
\affil 
ul.Miklukho-Maklaya 20-180, Moscow 117437 Russia\linebreak
E-mail: denis\@juriev.msk.ru
\endaffil
\date math.OC/9903001\enddate
\abstract\nofrills This note is devoted to an interactive game theoretic
formalization of dialogues as psycholinguistic phenomena and the unraveling 
of a hidden dialogue structure (the ``verbalization'') of 2-person 
differential interactive games. In the field-theoretic description of 
interactive games the dialogues are defined na{\"\i}vely as interactive games 
of discrete time with intention fields of continuous time; the correct 
mathematical formulation of such definition is proposed. The states and the 
controls of a dialogue correspond to the speech whereas the intention fields 
describe the understanding. In the case of dialogues the main inverse problem 
of interactive game theory is to describe geometrical and algebraical 
properties of the understanding. On the other hand, a precise mathematical 
definition of dialogues allows to formulate a problem of the unraveling of a 
hidden dialogue structure of any 2-person differential interactive game. Such 
procedure is called the verbalization. It means that the states of a
differential interactive game are interpreted as intention fields of a hidden 
dialogue and the problem is to describe such dialogue completely. If a 
2-person differential interactive game is verbalizable one is able to consider 
many linguistic (e.g. the formal grammar of a related hidden dialogue) or 
psycholinguistic (e.g. the dynamical correlation of various implications) 
aspects of it.
\endabstract
\endtopmatter
\document
The mathematical formalism of interactive games, which extends one of ordinary 
games [1] and is based on the concept of an interactive control, was proposed 
by the author [2] to take into account the complex composition of controls of 
a real human person, which are often complicated couplings of his/her 
cognitive and known controls with the unknown subconscious behavioral 
reactions. This formalism is applicable also to the description of external 
unknown influences and, thus, is useful for problems in computer science 
(e.g. the semi-artificially controlled distribution of resources) and 
mathematical economics (e.g. the financial games with unknown dynamical 
factors). This article is devoted to an embedding of dialogues as 
psycholinguistic phenomena into the interactive game theoretic picture. It 
will allow to apply the interactive game theoretic methods (such as {\sl a 
posteriori\/} analysis and short-term prediction procedures) to the dialogues 
and to adapt {\sl vice versa\/} the linguistic and psycholinguistic approaches 
to a wide class of the arbitrary 2-person interactive games of the nonverbal 
origin, especially to the computer interactive videogames. The linguistic and 
psycholinguistic methods may be crucial for the unraveling of the hiddenly 
verbal character of the accelerated naturally nonverbal cognitive computer 
and telecommunications based on droems and their dynamical reconstruction [4].

\head I. Interactive games\endhead

\subhead 1.1. Interactive systems and intention fields\endsubhead

\definition{Definition 1 [2]} {\it An interactive system\/} (with $n$
{\it interactive controls\/}) is a control system with $n$ independent controls 
coupled with unknown or incompletely known feedbacks (the feedbacks, which are
called the {\it behavioral reactions}, as well as their couplings with 
controls are of a so complicated nature that their can not be described 
completely). {\it An interactive game\/} is a game with interactive controls 
of each player.
\enddefinition

Below we shall consider only deterministic and differential interactive
systems. For symplicity we suppose that $n=2$. In this case the general
interactive system may be written in the form:
$$\dot\varphi=\Phi(\varphi,u_1,u_2),\tag1$$
where $\varphi$ characterizes the state of the system and $u_i$ are
the interactive controls:
$$u_i(t)=u_i(u_i^\circ(t),\left.[\varphi(\tau)]\right|_{\tau\leqslant t}),$$
i.e. the independent controls $u_i^\circ(t)$ coupled with the feedbacks on
$\left.[\varphi(\tau)]\right|_{\tau\leqslant t}$. One may suppose that the
feedbacks are integrodifferential on $t$.

\proclaim{Proposition [2]} Each interactive system (1) may be transformed
to the form (2) below (which is not, however, unique):
$$\dot\varphi=\tilde\Phi(\varphi,\xi),\tag2$$
where the magnitude $\xi$ (with infinite degrees of freedom as a rule) 
obeys the equation
$$\dot\xi=\Xi(\xi,\varphi,\tilde u_1,\tilde u_2),\tag3$$
where $\tilde u_i$ are the interactive controls of the form $\tilde 
u_i(t)=\tilde u_i(u_i^\circ(t); \varphi(t),\xi(t))$ (here the dependence
of $\tilde u_i$ on $\xi(t)$ and $\varphi(t)$ is differential on $t$, i.e.
the feedbacks are precisely of the form 
$\tilde u_i(t)=\tilde u_i(u_i^\circ(t);\varphi(t),\xi(t),
\dot\varphi(t),\dot\xi(t),\ddot\varphi(t),\ddot\xi(t),\ldots,
\varphi^{(k)}(t),\mathbreak\xi^{(k)}(t))$).
\endproclaim

\remark{Remark 1} One may exclude $\varphi(t)$ from the feedbacks in
the interactive controls $\tilde u_i(t)$. One may also exclude the
derivatives of $\xi$ and $\varphi$ on $t$ from the feedbacks.
\endremark

\definition{Definition 2 [2]} The magnitude $\xi$ with its dynamical equations
(3) and its cont\-ri\-bution into the interactive controls $\tilde u_i$ will 
be called {\it the intention field}.
\enddefinition

Note that the theorem holds true for the interactive games. In practice, the 
intention fields may be often considered as a field-theoretic description of 
subconscious individual and collective behavioral reactions. However, they 
may be used also the accounting of unknown or incompletely known external 
influences. Therefore, such approach is applicable to problems of computer 
science (e.g. semi-automatically controlled resource distribution) or 
mathematical economics (e.g. financial games with unknown factors).
The interactive games with the differential dependence of feedbacks are
called {\it differential}. Thus, the theorem states a possibility of
a reduction of any interactive game to a differential interactive game
by introduction of additional parameters -- {\sl the intention fields}.

\subhead 1.2. Differential interactive games and their 
$\varepsilon$--representation\endsubhead 

The most powerful way to investigate differential interactive games is
their {\sl a posteriori\/} analysis [3]. The $\varepsilon$--representation
of differential interactive games is a very convenient form of their
recording to perform such analysis. In the next paragraph 
$\varepsilon$--representation will be used to give the precise interactive
game theoretical definition of dialogues. 

\definition{Definition 3} The {\it $\varepsilon$--representation\/} of 
differential interactive game is a representation of the differential
feedbacks in the form
$$u_i(t)=u_i(u_i^\circ,\varphi(t),\ldots,\varphi^{(k)}(t);
\varepsilon_i(t))\tag4$$
with the {\sl known\/} function $u_i$ of all its arguments, where
the magnitudes $\varepsilon_i(t)\in\Cal E$ are {\sl unknown\/} functions of
$u_i^\circ$ and $\varphi(t)$ with its higher derivatives:
$$\varepsilon_i(t)=\varepsilon_i(u_i^\circ(t),\varphi(t),\dot\varphi(t),
\ldots,\varphi^{(k)}(t)).$$
\enddefinition

The short-term predictions based on {\sl a posteriori\/} analysis of
differential interactive games use some their $\varepsilon$--representations
with the frozen parameters $\varepsilon_i(t)$ at the moment $t_0$.

\head II. Dialogues and interactive games\endhead

\subhead 2.1. Dialogues as interactive games\endsubhead

Let us formalize dialogues as psycholinguistic phenomena in terms of
interactive games. First of all, note that one is able to consider
interactive games of discrete time as well as interactive games of
continuous time above.

\definition{Defintion 4A (the na{\"\i}ve definition of dialogues)}
The {\it dialogue\/} is a 2-person interactive game of discrete time with 
intention fields of continuous time.
\enddefinition

The states and the controls of a dialogue correspond to the speech whereas 
the intention fields describe the understanding. 

Let us give the formal mathematical definition of dialogues now.

\definition{Definition 4B (the formal definition of dialogues)}
The {\it dialogue\/} is a 2-person interactive game of discrete time of 
the form
$$\varphi_n=\Phi(\varphi_{n-1},\vec v_n,\xi(\tau)| 
t_{n-1}\!\leqslant\!\tau\!\leqslant\!t_n).\tag5$$
Here $\varphi_n\!=\!\varphi(t_n)$ are the states of the system at the
moments $t_n$ ($t_0\!<\!t_1\!<\!t_2\!<\!\ldots\!<\!t_n\!<\!\ldots$), 
$\vec v_n\!=\!\vec v(t_n)\!=\!(v_1(t_n),v_2(t_n))$ are the interactive 
controls at the same moments; $\xi(\tau)$ are the intention fields of 
continuous time with evolution equations
$$\dot\xi(t)=\Xi(\xi(t),\vec u(t)),\tag6$$
where $\vec u(t)=(u_1(t),u_2(t))$ are continuous interactive controls with 
$\varepsilon$--represented couplings of feedbacks:
$$u_i(t)=u_i(u_i^\circ(t),\xi(t);\varepsilon_i(t)).$$
The states $\varphi_n$ and the interactive controls $\vec v_n$ are certain
{\sl known\/} functions of the form
$$\aligned
\varphi_n=&\varphi_n(\vec\varepsilon(\tau),\xi(\tau)| 
t_{n-1}\!\leqslant\!\tau\!\leqslant\!t_n),\\
\vec v_n=&\vec v_n(\vec u^\circ(\tau),\xi(\tau)|
t_{n-1}\!\leqslant\!\tau\!\leqslant\!t_n).
\endaligned\tag7
$$
\enddefinition

Note that the most nontrivial part of mathematical formalization of dialogues
is the claim that the states of the dialogue (which describe a speech) are 
certain ``mean values'' of the $\varepsilon$--parameters of the intention
fields (which describe the understanding).

\remark{Remark 2} Note that in the case of dialogues the main inverse problem 
of interactive game theory [2] means to describe geometrical and algebraical 
properties of the understanding.
\endremark

The definition of dialogue may be generalized on arbitrary number of players.

\subhead 2.2. Unraveling a hidden dialogue structure of 2-person differential
interactive games. The verbalization\endsubhead 

An embedding of dialogues into the interactive game theoretical picture
generates the reciprocal problem: how to interpret an arbitrary differential
interactive game as a dialogue. Such interpretation will be called the
{\it verbalization}.

\definition{Definition 5} A differential interactive game of the form
$$\dot\varphi(t)=\Phi(\varphi(t),\vec u(t))$$
with $\varepsilon$--represented couplings of feedbacks 
$$u_i(t)=u_i(u^\circ_i(t),\varphi(t),\dot\varphi(t),\ddot\varphi(t),\ldots
\varphi^{(k)}(t);\varepsilon_i(t))$$
is called {\it verbalizable\/} if there exist {\sl a posteriori\/}
partition $t_0\!<\!t_1\!<\!t_2\!<\!\ldots\!<\!t_n\!<\!\ldots$ and the 
integrodifferential functionals
$$\aligned
\omega_n&(\vec\varepsilon(\tau),\varphi(\tau)|
t_{n-1}\!\leqslant\!\tau\!\leqslant\!t_n),\\
\vec v_n&(\vec u^\circ(\tau),\varphi(\tau)|
t_{n-1}\!\leqslant\!\tau\!\leqslant\!t_n)
\endaligned\tag8$$ 
such that
$$\omega_n=\Omega(\omega_{n-1},v_n;\varphi(\tau)|
t_{n-1}\!\leqslant\!\tau\!\leqslant\!t_n).
\tag 9$$
\enddefinition

The verbalizable differential interactive games realize a dialogue in sense
of Def.4.

The main heuristic hypothesis is that all differential interactive games
``which appear in practice'' are verbalizable. The verbalization means that 
the states of a differential interactive game are interpreted as intention 
fields of a hidden dialogue and the problem is to describe such dialogue 
completely. If a 2-person differential interactive game is verbalizable one 
is able to consider many linguistic (e.g. the formal grammar of a related 
hidden dialogue) or psycholinguistic (e.g. the dynamical correlation of 
various implications) aspects of it.

\remark{Remark 3} It will be very interesting to understand the possible
connections between quantization of interactive games [2] and their
verbalization.
\endremark

\remark{Remark 4} The verbalization may be an important part of the
strategical analysis of differential interactive games beyond the 
short-term predictions [3].
\endremark

\head III. Conclusions\endhead

Thus, the interactive game theoretic formalization of dialogues as 
psycholinguistic phenomena was performed and the unraveling of a hidden 
dialogue structure (the ``verbalization'') of 2-person differential 
interactive games was initiated. It will allow to apply the interactive game 
theoretic methods to the dialogues and to adapt {\sl vice versa\/} the 
linguistic and psycholinguistic approaches to a wide class of the arbitrary 
2-person interactive games of the nonverbal origin.

\Refs
\roster
\item"[1]" Isaaks R., Differential games. Wiley, New York, 1965;\newline
Owen G., Game theory, Saunders, Philadelphia, 1968.
\item"[2]" Juriev D., Interactive games and representation theory. I,II.
E-prints: math.FA/9803020, math.RT/9808098.
\item"[3]" Juriev D., The laced interactive games and their {\sl a
posteriori\/} analysis. E-print:\linebreak math.OC/9901043; Differential
interactive games: The short-term predictions. E-print: math.OC/9901074.
\item" [4]" Juriev D., Droems: experimental mathematics, informatics and
infinite dimensional geometry. Report RCMPI-96/05+ [e-version: 
cs.HC/9809119].
\endroster
\endRefs
\enddocument